\documentclass[10pt]{amsart}
\usepackage{amsfonts,amssymb,amsthm,eucal,amsmath}
\allowdisplaybreaks
\newtheorem{thm}{Theorem}

\newcommand{\R}{\mathbb{R}}

\newcommand{\E}{\mathbb{E}}

\newcommand{\inprod}[2]{\left\langle #1, #2 \right\rangle}

\newcommand{\ds}{\displaystyle}

\renewcommand{\P}{\mathbb{P}}

\renewcommand{\L}{\mathcal{L}}

\begin{document}

\title{Two multivariate central limit theorems}
\author{Elizabeth Meckes}
\address{ Elizabeth Meckes\\
American Institute of Mathematics and
Department of Mathematics, Cornell University\\
Ithaca, NY 14853-4201.}
\email{esmeckes@math.cornell.edu}

\begin{abstract}
In this paper, explicit error bounds are derived in the approximation
of rank $k$ projections of certain $n$-dimensional random vectors by 
standard $k$-dimensional Gaussian random vectors.  The bounds are given
in terms of $k$, $n$, and a basis of the $k$-dimensional space onto 
which we project.  The random vectors considered are two generalizations
of the case of a vector with independent, identically distributed
 components.  In the first case, the 
random vector has components which are independent but need not have the
same distribution.  The second case deals with finite exchangeable sequences
of random variables.
\end{abstract}
\maketitle

\section{Introduction}
The classical central limit theorem says that, under mild
conditions, the random variable
$S_n=\frac{1}{\sqrt{n}}\sum_{i=1}^nX_i$ is approximately Gaussian, for 
a sequence of $n$ independent, identically distributed random variables 
$X_i$ and $n$ large.  That is, if $X$ is a random vector of $\R^n$ with i.i.d. components,
then the orthogonal projection of $X$ in the direction $(1,1,\ldots,1)$
is approximately Gaussian.  It is natural to ask for what other directions
the projection of such a random vector is approximately Gaussian.  In 
particular, the Berry-Esseen theorem for sums of independent, non-identically
distributed random variables implies that for $X$ 
as above, $$\sup_{t\in\R}\big| \P[\inprod{\theta}{X} \le t] - \Phi(t) \big|
  \le 0.8 \Big(\E |X_i|^3 \Big) 
  \sum_{j=1}^n |\theta_j|^3,$$
where $\Phi$ denotes the standard normal distribution function and $\theta$
is any unit vector in $\R^n$.  Thus $\inprod{\theta}{X}$ is close to 
Gaussian as long as $ \sum_{j=1}^n |\theta_j|^3$ is small.  Roughly, this
happens as long as there aren't a small number of coordinates of $\theta$
controlling the value of $\inprod{\theta}{X}$; i.e., the components of 
$\theta$ are all of similar size.  More generally, one could ask
when higher rank projections of $X$ are close to Gaussian; that is, not
only consider the asymptotic normality of individual projections of $X$, 
but the asymptotic independence of projections in different directions.  
The basic result of this paper is the following quantitative bound
on the distance from a rank $k$ projection of $X$ to a standard Gaussian 
random vector in a fixed dimension, where distance is measured here by 
comparing the integrals of $C^2$ test functions.  In the following theorem,
$C^2_c(\R^k)$ denotes the space of compactly supported, real-valued functions
on $\R^k$ with two continuous derivatives; $\||\nabla g(x)|\,\|_\infty$ is the
maximum length of the gradient of $g$ and
$\ds |g|_2=\max_{i,j}\sup_x\left|\frac{\partial^2 g}{\partial x_i\partial x_j}
(x)\right|$.  
The $\ell_p$ norm of a vector $\theta\in\R^n$ is denoted $\|\theta\|_p=
\left(\sum_{i=1}^n|\theta_i|^p\right)^{\frac{1}{p}}.$
\begin{thm}\label{simple}
Let $X_1,\ldots,X_n$ be independent, identically distributed random
variables with $\E X_1=0$ and $\E X_1^2=1$.
Let $\theta_1,\ldots,\theta_k$
be fixed vectors in $\R^n$ with $\theta_i=(\theta_i^1,\ldots,\theta_i^n)$, such that
$\inprod{\theta_i}{\theta_j}=\delta_{ij}.$
Define a random vector $S_n\in\R^k$ by
$$S_n^i=\inprod{\theta_i}{X}=\sum_{r=1}^n\theta_i^rX_r.$$
Then for $g\in C_c^2(\R^k),$
\begin{equation}\begin{split}
\big|\E g(S_n)-\E g(Z)\big|&\le\frac{\sqrt{k}}{2}\||
\nabla g|\|_\infty\sqrt{\E X_1^4-1}
\left(\sum_{i=1}^k\|\theta_i\|_4^2\right)\\&\qquad+\frac{4}{3}k^2|g|_2
\left(\E|X_1|^3\right)\left(\sum_{i=1}^k
\|\theta_i\|_3^3\right),\end{split}\end{equation}
where $Z$ is distributed as a standard normal random vector in $\R^k$.
\end{thm}
The following example is useful to consider.  Suppose the $\theta_i$ are such 
that $|\theta_i^j|=\frac{1}{\sqrt{n}}$
for each $j$; i.e., the $\theta_i$ are orthogonal 
unit vectors in the directions of
corners of the hypercube.  As long as $n$ is large and a multiple of 4,
there are more such vectors than we can make use of -- see \cite{dL}.
Then the norm expressions in the bound above 
reduce to
\begin{equation*}
\sum_{i=1}^k\|\theta_i\|_4^2=\frac{k}{\sqrt{n}},\phantom{\sum_{i=1}^k
\|\theta_i\|_4^2=}
\sum_{i=1}^k
\|\theta_i\|_3^3=\frac{k}{\sqrt{n}}.
\end{equation*}
Thus for directions chosen in this way, projections of rank $k$ are close 
to Gaussian as long as $k=o\left(n^{1/6}\right)$.

Furthermore, if $\theta_i$ is random on the sphere, then 
$$\E \|\theta_i\|_4^2\le\sqrt{\E\sum_{r=1}^n(\theta_i^r)^4}=\sqrt{\frac{3}{
n+2}}$$
and $$\E\|\theta_i\|_3^3=\sum_{r=1}^n\E|\theta_i^r|^3=\frac{n\Gamma\left(\frac{
n}{2}\right)}{\sqrt{\pi}\Gamma\left(\frac{n}{2}+\frac{3}{2}\right)}\approx
\sqrt{\frac{8}{n\pi}}.$$
(See \cite{foll} for a straightforward approach to integrating even-degree
monomials over the sphere; the proof given there extends to odd-degree
monomials in the absolute values of coordinates as well.)
It follows that for the $\theta_i$ chosen at random (subject to the 
orthogonality condition), there are absolute constants $c_1$ and $c_2$  
such that for every $g\in C_c^2(\R^k)$,
$$\E_{\bf \theta}\left(\big|\E g(S_{n,{\bf \theta}})
-\E g(Z)\big|
\right)\le
\frac{1}{\sqrt{n}}\Big[c_1k^{3/2}\||\nabla g|\|_\infty\sqrt{\E X_1^4-1}+
c_2k^3|g|_2\E|X_1|^3\Big].$$
This implies that a typical projection of rank $k$ 
is close to Gaussian for $k=o(n^{1/6})$.

\medskip

Theorem \ref{simple} is generalized below in two directions.  
In the following version, the $X_i$ are assumed to be independent, but
need not be identically distributed.
\begin{thm}\label{orth}
Let $X_1,\ldots,X_n$ be independent (not necessarily identically 
distributed) random
variables with $\E X_i=0$ and $\E X_i^2=1$ for each $i$.  
Let $\theta_1,\ldots,\theta_k$
be fixed vectors in $\R^n$ with $\theta_i=(\theta_i^1,\ldots,\theta_i^n)$, such that
$\inprod{\theta_i}{\theta_j}=\delta_{ij}.$
Define a random vector $S_n\in\R^k$ by
$$S_n^i=\inprod{\theta_i}{X}=\sum_{r=1}^n\theta_i^rX_r.$$
Then for $g\in C_c^2(\R^k),$
\begin{equation}\begin{split}\label{spbd}
\big|\E g(S_n)-\E g(Z)\big|&\le\frac{\sqrt{k}}{2}\||
\nabla g|\|_\infty\sqrt{\max_{1\le i\le n}\E X_1^4-1}
\left(\sum_{i=1}^k\|\theta_i\|_4^2\right)\\&\qquad+\frac{4}{3}k^2|g|_2
\left(\max_{1\le j\le n}\E|X_j|^3\right)\left(\sum_{i=1}^k
\|\theta_i\|_3^3\right),\end{split}\end{equation}
where $Z$ is distributed as a standard normal random vector in $\R^k$.
\end{thm}

\medskip

\medskip

Theorem \ref{orth} can be generalized further to require the vectors $\theta_i$
only to be linearly independent.  Let $H_g(x)$ be the Hessian matrix of $g$
at $x$, and let $$\|\|H_g\|_{op}\|_\infty=\sup_x\|H_g(x)\|_{op}$$
where $\|A\|_{op}$ is the operator norm of the matrix $A$.  Thus $
\|\|H_g\|_{op}\|_\infty$ is the supremum over $x$ of the largest 
eigenvalue (in absolute value) of $H_g(x)$.
\begin{thm}\label{non-orth}
Let $X_1,\ldots,X_n$ be independent (not necessarily identically
distributed) random variables with 
$\E X_i=0$ and $\E X_i^2=1$ for each $i$.  Let $\theta_1,\ldots,\theta_k$
be fixed, linearly independent vectors in $\R^n$
with $\theta_i=(\theta_i^1,\ldots,
\theta_i^n)$, such that
$\|\theta_i\|_2=1$ for each $i$.  Let $c_{ij}=\inprod{
\theta_i}{\theta_j}.$  Define a random vector $S_n\in\R^k$ by
$$S_n^i=\inprod{\theta_i}{X}=\sum_{r=1}^n\theta_i^rX_r$$
and let $\widetilde{Z}$ be a Gaussian random vector with covariance matrix
$C=(c_{ij})_{i,j=1}^k.$  Then for $f\in C^2_c(\R^k),$
\begin{equation}\begin{split}\label{genbd}
\big|\E f(S_n)-\E f(\widetilde{Z})\big|&\le\frac{1}{2}\sqrt{
\lambda k}\||
\nabla f|\|_\infty\sqrt{\max_{1\le i\le n}\E X_i^4-1}
\left(\sum_{i=1}^k\|\theta_i\|_4^2\right)\\&\qquad
+\frac{4}{3}\lambda k^2\Big(\|\|H_f\|_{op}\|_\infty\Big)
\left(\max_{1\le i\le n}\E|X_i|^3\right)\left(\sum_{i=1}^k
\|\theta_i\|_3^3\right),\end{split}\end{equation}
where $\lambda$ is the largest eigenvalue of $C.$
\end{thm}

Theorem \ref{non-orth} follows from Theorem \ref{orth} using 
a fairly straightforward linear algebra argument.

\medskip

Theorem \ref{simple} can also be generalized in a different direction, by
weakening the independence assumption.  In the following version, the 
sequence $X_1,\ldots,X_n$ is assumed to be exchangeable, i.e., 
$(X_1,\ldots,X_n)\stackrel{\mathcal{L}}{=}(X_{\sigma(1)}\ldots,X_{\sigma(n)})$
for any permutation $\sigma$, but $(X_1,\ldots,X_n)$ need not
have independent entries.  Theorem \ref{exch-thm} is not a generalization
of Theorem \ref{simple} in the strictest sense, as it has the additional
technical requirement that $\sum_r\theta_i^r=0$ for each $i$.  In 
what follows, let $|g|_1=\max_{1\le i\le k}\left\|\frac{\partial g}{\partial 
x_i}\right\|_{\infty}.$
\begin{thm}\label{exch-thm}
Let $(X_1,\ldots,X_n)$ be a finite exchangeable sequence of random variables
with $\E X_1=0$ and $\E X_1^2=1$.  Let $\{\theta_i\}_{i=1}^k$ be an 
orthonormal set of vectors in $\R^n$, such that $\sum_{r=1}^n\theta_i^r=0$
for each $i$.  Define the random vector $S_n$ in
$\R^k$ by 
$$S_n^i=\inprod{\theta_i}{X}=\sum_{r=1}^n\theta_i^rX_r.$$
Then there are absolute constants $a$, $b$, $c$ such that for any 
$g\in C^2_c(\R^k),$
\begin{equation}\begin{split}
\big|\E g(S_n)-\E g(Z)\big|&\le
ak|g|_1\left(\sqrt{\big|\E X_1X_2X_3X_4\big|}+\sqrt{\big|\E(X_1^2-1)(X_2^2-
1)\big|}\right)\\&\quad
+b|g|_1\sqrt{\E X_1^4}\left(\sum_{i=1}^k\|\theta_i\|_4\right)^2+
ck^2|g|_2\E|X_1|^3\left(\sum_{i=1}^k\|\theta_i\|_3^3\right).\end{split}
\end{equation}
\end{thm}

In the case that the entries are independent, the first two error
terms vanish; one can interpret their presence as a requirement that the
dependence among the $X_i$ must be weak.

In the same way as one obtains Theorem \ref{non-orth} from Theorem \ref{orth},
one can weaken the orthonormality requirement on the $\theta_i$ of 
Theorem \ref{exch-thm} to the requirement that they be linearly independent.
This yields the following.
\begin{thm}\label{exch-lin-ind}
Let $(X_1,\ldots,X_n)$ be an exchangeable sequence of random variables
with $\E X_1=0$ and $\E X_1^2=1$.  Let $\{\theta_i\}_{i=1}^k$ be a linearly independent set of vectors in $\R^n$, 
such that $\sum_{r=1}^n\theta_i^r=0$
for each $i$.  Let $c_{ij}=\inprod{\theta_i}{\theta_j}.$   
Define the random vector $S_n$ in
$\R^k$ by 
$$S_n^i=\inprod{\theta_i}{X}=\sum_{r=1}^n\theta_i^rX_r,$$
and let $\widetilde{Z}$ be a Gaussian random vector with covariance matrix
$C=(c_{ij})_{i,j=1}^k.$
Then there are absolute constants $a$, $b$, $c$ such that for any 
$g\in C^2_c(\R^k),$
\begin{equation}\begin{split}
\big|\E g(S_n)-\E g(\widetilde{Z})\big|&\le
ak\sqrt{\lambda}\||\nabla g|\|_\infty\left(\sqrt{\big|\E X_1X_2X_3X_4\big|}+
\sqrt{\big|\E(X_1^2-1)(X_2^2-
1)\big|}\right)\\&\quad
+b\sqrt{\lambda}\||\nabla g|\|_\infty\sqrt{\E X_1^4}\left(\sum_{i=1}^k\|
\theta_i\|_4\right)^2+
ck^2\lambda\Big(\|\|H_g\|_{op}\|_\infty\Big)\E|X_1|^3\left(\sum_{i=1}^k\|\theta_i\|_3^3\right),\end{split}
\end{equation}
where $\lambda$ is the largest eigenvalue of $C$.
\end{thm}

\section{Proofs}

\begin{proof}[Proof of Theorem \ref{non-orth} from Theorem \ref{orth}]
Perform the Gram-Schmidt algorithm on the set of vectors $\{\theta_i\}$: 
since the $\theta_i$ are linearly independent, there is an invertible matrix
$B$ such that for $\eta_i:=\sum_jB^{-1}_{ij}\theta_j,$  $\inprod{\eta_i}{
\eta_j}=\delta_{ij}$.
By assumption,
\begin{equation*}\begin{split}
c_{ij}&=\inprod{\theta_i}{\theta_j}\\&=\inprod{\sum_pB_{ip}\eta_p}{\sum_qB_{jq}
\eta_q}\\&=\sum_{p,q}B_{ip}B_{jq}\inprod{\eta_p}{\eta_q}\\&=\sum_pB_{ip}
B_{jp}.\end{split}\end{equation*}
Thus $BB^T=C.$

Now, let $f:\R^k\to\R$ and define $h:\R^k\to\R$ by $h(x)=f(Bx).$  Define
$\widetilde{S}_n^i:=\sum_{r=1}^n\eta_i^rX_r$.  Then
$$(B\widetilde{S}_n)_i=\sum_jB_{ij}\inprod{\eta_j}{X}=
\inprod{\sum_jB_{ij}\eta_j}{X}=\inprod{\theta_i}{X}=S_n^i,$$
and so $h(\widetilde{S}_n)=f(B\widetilde{S}_n)=f(S_n).$
If $Z$ is a standard Gaussian random vector, then $h(Z)=f(BZ)$ and
$BZ$ is a Gaussian random vector with covariance matrix $BB^T=C.$
Applying Theorem \ref{orth} for this test function $h$ to the random
vector $\widetilde{S}_n^i$ yields
\begin{equation*}\begin{split}
\big|\E f(S_n)-\E f(BZ)\big|&=
\big|\E h(\widetilde{S}_n)-\E h(Z)\big|\\&\le \frac{1}{2}k^{1/2}\||
\nabla h|\|_\infty\sqrt{\max_{1\le j\le n}\E X_j^4-1}\left(\sum_{i=1}^k\|
\theta_i\|_4^2\right)\\&\qquad
+\frac{4}{3}k^2|h|_2\left(\max_{1\le j\le n}\E|X_j|^3\right)\left(\sum_{i=1}^k
\|\theta_i\|_3^3\right),\end{split}\end{equation*}
so it remains to estimate $\||\nabla h|\|_\infty$ and $|h|_2.$

To estimate $\||\nabla h|\|_\infty,$ first note that if $B$ is viewed as
an operator on $\R^k$, then its operator norm is its largest singular value.  
That is, $$\|B x\|\le\sqrt{\lambda}
\|x\|,$$ where $\lambda$ is the largest
eigenvalue of $C.$  Since $\||\nabla h|\|_\infty$ is the 
Lipschitz constant of $h$, it follows that $\||\nabla h|\|_\infty\le\sqrt{
\lambda}\||\nabla f|\|_\infty.$

To estimate $|h|_2$, note that
$$\frac{\partial^2h}{\partial x_\ell\partial x_p}=\sum_{r,s}B_{rp}B_{s\ell}
\frac{\partial^2f}{\partial x_s\partial x_r}(Bx)=(B^TH_f(Bx)B)_{\ell p},$$
where $(H_f)_{ij}=\left(\frac{\partial^2 f}{\partial x_i\partial x_j}
\right)_{ij}$ is the Hessian of $f$.
Now,
$$|B^THB|_{\ell p}=|\inprod{B^THBe_\ell}{e_p}|=|\inprod{HBe_\ell}{Be_p}|
\le \|H\|_{op}\|B\|_{op}^2.$$
As stated above, $\|B\|_{op}^2\le\lambda,$ the largest eigenvalue of $C;$
this completes the proof.

\end{proof}

\medskip

{\bf Remarks:}
\begin{enumerate}
\item  To obtain a bound in Theorem \ref{non-orth} which doesn't 
involve the operator norm of the Hessian of $f$, one can 
estimate the operator norm by the Hilbert-Schmidt norm:
$$\|H_f(Bx)\|_{op}\le\sqrt{\sum_{r,s=1}^k\left(\frac{\partial^2 f}{\partial
x_r\partial x_s}(Bx)\right)^2}\le k|f|_2.$$
\item The proof of Theorem \ref{exch-lin-ind} from Theorem \ref{exch-thm}
is exactly the same as the proof above.
\end{enumerate}
\medskip

The proofs of Theorems \ref{orth} and \ref{exch-thm} are applications of
the following abstract normal approximation theorem, proved in \cite{cm}
\begin{thm}\label{discrete}
Let $X$ and $X'$ be two random vectors in $\R^k$ such that $\L(X)=\L(X')$,
and let $Z\in\R^k$ be a standard Gaussian random vector.
Suppose there is a constant $\lambda$ and random variables $E_{ij}$ such that
\begin{enumerate}
\item $\E\left[X_i'-X_i\big|X\right]=-\lambda X_i$\label{lin-cond}
\item $\E\left[(X_i'-X_i)(X_j'-X_j)\big|X\right]=2\lambda\delta_{ij}+
E_{ij}.$\label{quad-cond}
\end{enumerate}
Then if $g\in C^2(\R^k),$
\begin{equation}\begin{split}\label{discbound}
\big|\E g(X) -\E g(Z)\big| &\le\min\left\{\frac{|g|_1}{2\lambda}\sum_{i,j}
\E|E_{ij}|, \frac{\sqrt{k}\||\nabla g|\|_\infty}{2\lambda}\E\biggl(\sum_{i,j}
E_{ij}^2\biggr)^{1/2}\right\}\\
&\quad +
\frac{k^2|g|_2}{6\lambda}\sum_i\E|X_i'-X_i|^3.
\end{split}\end{equation}

\end{thm}

In the contexts in which Theorem \ref{discrete} is applied below, the 
pair of vectors $X,X'$ will in fact be constructed not only 
to have the same law but to be exchangeable.

\begin{proof}[Proof of Theorem \ref{orth}]
From the random vector $S_n$ with $S_n^i=\inprod{\theta_i}{X}$, make an exchangeable pair of vectors
$(S_n,S_n')$ by choosing $I\in\{1,\ldots,n\}$ at random,
independent of $\{X_i\}$, and replacing $X_I$ by an independent copy
$X_I^*.$  That is,
$$(S_n')^i=S_n^i-\theta_i^IX_I+\theta_i^IX_I^*.$$
Then
\begin{equation*}
\E\left[(S_n')^i-S_n^i\big|\{X_j\}_{j=1}^n\right]=-\frac{1}{n}
\sum_{j=1}^n\theta_i^jX_j=-\frac{1}{n}S_n^i,
\end{equation*}
thus the proportionality condition of Theorem \ref{discrete} holds with
$\lambda=\frac{1}{n}$.

Next,
\begin{equation*}\begin{split}
\E\left[\left((S_n')^i-S_n^i\right)^2\big|\{X_j\}_{j=1}^n\right]&=\E\left[
\left.(\theta_i^I)^2(X_I-X_I^*)^2\right|\{X_j\}_{j=1}^n\right]\\
&=\frac{1}{n}\sum_{\ell=1}^n\E\left[(\theta_i^\ell)^2(X_\ell-X_\ell^*)^2
\big|\{X_j\}_{j=1}^n\right]
\\&=\frac{1}{n}\sum_{\ell=1}^n(\theta_i^\ell)^2(X_\ell^2+1)\\
&=\frac{2}{n}+\frac{1}{n}\sum_{\ell=1}^n(\theta_i^\ell)^2(X_\ell^2-1),
\end{split}\end{equation*}
since $\sum_{\ell}(\theta_i^\ell)^2=1.$  Thus one can take
$E_{ii}=\frac{1}{n}\sum_{\ell=1}^n(\theta_i^\ell)^2(X_\ell^2-1)$.

If $i\neq j$,
\begin{equation*}\begin{split}
\E\left[\left((S_n')^i-S_n^i\right)\left((S_n')^j-S_n^j\right)\big|
\{X_l\}_{l=1}^n\right]&=\E\left[\left.\theta_i^I\theta_j^I
(X_I-X_I^*)^2\right|\{X_l\}_{l=1}^n\right]\\&=\frac{1}{n}\sum_{r=1}^n
\theta_i^r\theta_j^r\E\left[(X_r-X_r^*)^2\big|\{X_l\}_{l=1}^n\right]\\
&=\frac{1}{n}\sum_{r=1}^n\theta_i^r\theta_j^r(X_r^2+1)\\
&=\frac{1}{n}\sum_{r=1}^n\theta_i^r\theta_j^r(X_r^2-1),
\end{split}\end{equation*}
where the last line follows because $\sum_r\theta_i^r\theta_j^r=0.$
Thus $E_{ij}=\frac{1}{n}\sum_{r=1}^n\theta_i^r\theta_j^r(X_r^2-1)$
for all $i$ and $j$.
Now,
\begin{equation*}\begin{split}
\sum_{i,j=1}^kE_{ij}^2&=\frac{1}{n^2}\sum_{i,j=i}^k\sum_{\ell,r=1}^n
\theta_i^\ell\theta_i^r\theta_j^\ell\theta_j^r(X_\ell^2-1)(X_r^2-1).
\end{split}\end{equation*}
Making use of the facts that $\E (X_\ell^2-1)=0$ and $X_\ell$ and $X_r$ are
independent for $\ell\neq r$, together with the Cauchy-Schwarz inequality
yields
\begin{equation*}\begin{split}
\E\sqrt{\sum_{i,j=1}^kE_{ij}^2}&=\frac{1}{n}\E\sqrt{\sum_{i,j=i}^k
\sum_{\ell,r=1}^n\theta_i^\ell\theta_i^r\theta_j^\ell\theta_j^r(X_\ell^2-1)(X_r^2-1)}\\
&\le\frac{1}{n}\sqrt{\E\left(\sum_{i,j=i}^k\sum_{\ell,r=1}^n
\theta_i^\ell\theta_i^r\theta_j^\ell\theta_j^r(X_\ell^2-1)(X_r^2-1)\right)}\\
&=\frac{1}{n}\sqrt{\E\left(\sum_{i,j=1}^k\sum_{\ell=1}^n
(\theta_i^\ell)^2(\theta_j^\ell)^2(X_\ell^2-1)^2\right)}\\
&\le\frac{1}{n}\sqrt{\max_{1\le i\le n}(\E X_i^4-1)\left(\sum_{i,j=1}^k\sum_{\ell=1}^n
(\theta_i^\ell)^2(\theta_j^\ell)^2\right)}.
\end{split}\end{equation*}
Now,
\begin{equation*}
\sum_{\ell=1}^n(\theta_i^\ell)^2(\theta_j^\ell)^2\le\sqrt{\sum_{\ell=1}
^n(\theta_i^\ell)^4}\sqrt{\sum_{\ell=1}^n(\theta_j^\ell)^4}=
\|\theta_i\|_4^2\|\theta_j\|_4^2,
\end{equation*}
and so
$$\E\sqrt{\sum_{i,j=1}^kE_{ij}^2}\le\frac{1}{n}\left(\sum_{i=1}^k\|\theta_i
\|_4^2\right)\sqrt{\max_{1\le i\le n}\E X_i^4-1}.$$
Finally,
\begin{equation*}\begin{split}
\E\big|(S_n')^i-S_n^i\big|^3&=\frac{1}{n}\sum_{j=1}^n|\theta_i^j|^3
\E|X_j-X_j^*|^3\\
&\le\frac{8}{n}\sum_{j=1}^n|\theta_i^j|^3\E|X_j|^3\\&\le\frac{8\|\theta_i
\|_3^3}{n}\max_{1\le i\le n}\E|X_i|^3,
\end{split}\end{equation*}
where the second line follows from the $L_3$ triangle inequality and the fact
that $X_j^*$ has the same distribution as $X_j$.  The statement of
the theorem is now an immediate consequence of Theorem \ref{discrete}.
\end{proof}

\bigskip

\begin{proof}[Proof of Theorem \ref{exch-thm}]
Starting from $S_n$, make an exchangeable pair of random vectors as follows.  Choose a pair
of indices $I\neq J$ at random from $\{1,\ldots,n\}.$  Let $\tau=(IJ),$
the permutation on $n$ letters that transposes $I$ and $J$, and let
$$X'=(X_{\tau(1)},\ldots,X_{\tau(n)}).$$
Then $$(S_n')^i=\inprod{\theta_i}{X'}=S_n^i+\theta_i^IX_J-\theta_i^IX_I+
\theta_i^JX_I-\theta_i^JX_J,$$
that is,
\begin{equation}\label{diff}
(S_n')^i-S_n^i=(\theta_i^I-\theta_i^J)(X_J-X_I).\end{equation}
Let $\sideset{}{'}\sum$ denote summing over distinct indices.  Then
\begin{equation*}\begin{split}
\E\left[(S_n')^i-S_n^i\big|X\right]&=\frac{1}{n(n-1)}\sideset{}{'}\sum_{r,s}
\big(\theta_i^rX_s-\theta_i^rX_r+\theta_i^sX_r-\theta_i^sX_s\big)\\
&=-\frac{2}{n-1}\sum_r\theta_i^rX_r\\&=-\frac{2}{n-1}S_n^i,
\end{split}\end{equation*}
where the second line follows as $\sum_r\theta_i^r=0.$
Thus the proportionality condition of Theorem \ref{discrete} holds with $\lambda =\frac{2}{n-1}$.

The next step is to compute and bound the error terms $E_{ij}.$  First,
consider $i=j$.  From \eqref{diff},
\begin{equation*}\begin{split}
\E\left[((S_n')^i-S_n^i)^2\big|X\right]&=\frac{1}{n(n-1)}
\sideset{}{'}\sum_{r,s}(\theta_i^r-\theta_i^s)^2(X_s-X_r)^2\\&=
\frac{1}{n(n-1)}\sum_{r,s}(\theta_i^r-\theta_i^s)^2(X_s-X_r)^2\\&=
\frac{1}{n(n-1)}\left[2\sum_{r,s}(\theta_i^r)^2X_s^2-4\sum_{r,s}(\theta_i^r)^2
X_rX_s+2n\sum_r(\theta_i^r)^2X_r^2\right.\\&\qquad\qquad\qquad\qquad\left.
-4\sum_{r,s}\theta_i^r\theta_i^sX_s^2+4
\sum_{r,s}\theta_i^rX_r\theta_i^sX_s\right]\\&=\frac{2}{n(n-1)}
\left[\sum_rX_r^2+n\sum_r(\theta_i^r)^2X_r^2-2\left(\sum_r(\theta_i^r)^2X_r
\right)\left(\sum_sX_s\right)+2(S_n^i)^2\right]\\&=\frac{4}{n-1}+
\frac{2}{n(n-1)}\left[\left(\sum_rX_r^2-n\right)+n\left(\sum_r(\theta_i^r)^2
X_r^2-1\right)\right.\\&\qquad\qquad\qquad\qquad\left.-2\left(\sum_r(
\theta_i^r)^2X_r\right)\left(\sum_sX_s\right)+2(S_n^i)^2\right].
\end{split}\end{equation*}
The error $E_{ii}$ can thus be taken to be
\begin{equation*}
E_{ii}=\frac{2}{n(n-1)}\left[\sum_r(X_r^2-1)+n\sum_r(\theta_i^r)^2
(X_r^2-1)-2\left(\sum_r(
\theta_i^r)^2X_r\right)\left(\sum_sX_s\right)+2(S_n^i)^2\right].
\end{equation*}
To bound $\E|E_{ii}|,$ first apply the triangle inequality and treat each of
the four terms above separately.  First,
\begin{equation*}\begin{split}
\E\left|\sum_r(X_r^2-1)\right|&\le\sqrt{\E\left(\sum_{r,s}(X_r^2-1)(X_s^2-1)
\right)}\\&\le\sqrt{n\big|\E\left(X_1^4-1\right)\big|}+\sqrt{n(n-1)\big|
\E\left[(X_1^2-1)(X_2^2-1)\right]\big|}.
\end{split}\end{equation*}
Next,
\begin{equation*}\begin{split}
\E\left|\sum_r(\theta_i^r)^2(X_r^2-1)\right|&\le\sqrt{\E\left(\sum_{r,s}
(\theta_i^r)^2(\theta_i^s)^2(X_r^2-1)(X_s^2-1)\right)}\\&\le\sqrt{\left|
\E\left[X_1^4-1\right]\right|\sum_r(\theta_i^r)^4+
\left|\E\left[(X_1^2-1)(X_2^2-1)\right]\right|\sum_{r,s}
(\theta_i^r)^2(\theta_i^s)^2}\\&\le\|\theta_i\|_4^2\sqrt{\big|\E X_1^4-1\big|}+
\sqrt{\big|\E\left[(X_1^2-1)(X_2^2-1)\right]\big|}.
\end{split}\end{equation*}
Note also that the normalization is such that $\E (S_n^i)^2=1,$ thus only 
the second-last term remains to be estimated.  
As before, start by applying the Cauchy-Schwarz inequality:
$$\E\left|\sum_{r,s}(\theta_i^r)^2X_rX_s\right|\le\sqrt{\E\sum_{k,\ell,r,s}
(\theta_i^k)^2(\theta_i^r)^2X_kX_\ell X_rX_s}.$$
Now, breaking up the sum by the equality structure of the indices and using
the exchangeability gives
\begin{equation*}\begin{split}
&\sqrt{\E\sum_{k,\ell,r,s}(\theta_i^k)^2(\theta_i^r)^2X_kX_\ell X_rX_s}\\
&\qquad\le\Big(n^2\big|\E\big[X_1X_2X_3X_4\big]\big|+n^2\|\theta_i\|_4^4
\big|\E\big[X_1^2X_2X_3\big]\big|+n\|\theta_i\|_4^4\E\big[X_1^2X_2^2\big]+
2n\|\theta_i\|_4^4\big|\E\big[X_1^3X_2\big]\big|\\&\qquad\qquad
+\|\theta_i\|_4^4\E\big[X_1^4\big]+n\big|\E\big[X_1^2X_2X_3\big]\big|+
2\big|\E\big[X_1^3X_2\big]\big|+2n\big|\E\big[X_1^2X_2X_3\big]\big|+2
\E\big[X_1^2X_2^2\big]
\Big)^{\frac{1}{2}}\\&\qquad\le n\left[\sqrt{\big|\E X_1X_2X_3X_4\big|}+
\|\theta_i\|_4^2\sqrt{\big|\E X_1^2X_2X_3\big|}\right]\\&\qquad\qquad+
\sqrt{n}\left[\|\theta_i\|_4^2\left(\sqrt{2\big|\E X_1^3X_2\big|}+\sqrt{\E
X_1^2X_2^2}\right)+
\sqrt{3\big|\E X_1^2X_2X_3\big|}\right]\\&\qquad\qquad\qquad
+\left[\|\theta_i\|_4^2\sqrt{\E X_1^4}
+\sqrt{2\big|\E X_1^3X_2\big|}+\sqrt{2\E X_1^2X_2^2}\right].
\end{split}\end{equation*}
By H\"older's inequality and the 
exchangeability of the sequence, 
$$\big|\E X_1^3X_2\big|\le \E X_1^4,\qquad\big|\E X_1^2 X_2^2\big|\le\E X_1^4,
\qquad{\rm and }\,\big|\E X_1^2X_2X_3\big|\le\E X_1^4.$$
Also,
$$\frac{1}{\sqrt{n}}\le\|\theta_i\|_4^2\le 1$$
since $\|\theta_i\|_2=1.$  Thus there are constants $c,c'$ such that 
$$\sqrt{\E\sum_{k,\ell,r,s}(\theta_i^k)^2(\theta_i^r)^2X_kX_\ell X_rX_s}
\le cn\sqrt{\big|\E X_1X_2X_3X_4\big|}+c'n\|\theta_i\|_4^2\sqrt{
\E X_1^4}.$$
All together, this shows that there are constants $c_1,c_2,c_3$ such that 
\begin{equation}\label{Eii}
\E|E_{ii}|\le \frac{c_1}{n}\sqrt{\big|\E X_1X_2X_3X_4\big|}+\frac{c_2}{n}
\sqrt{\big|\E\left(X_1^2-1\right)\left(X_2^2-1\right)\big|}+\frac{c_3}{n}
\|\theta_i\|_4^2\sqrt{\E X_1^4}.
\end{equation}

\bigskip

Next, consider $E_{ij}$ for $i\neq j$.  From \eqref{diff},
\begin{equation*}\begin{split}
\E\left[\left.\left((S_n')^i-S_n^i\right)\left((S_n')^j-S_n^j\right)\right|
\{X_\ell\}_{\ell=1}^n\right]&=\E\left[\left.\left(\theta_i^I-
\theta_i^J\right)\left(\theta_j^I-\theta_j^J\right)\left(X_J-X_I\right)^2
\right|\{X_\ell\}\right]\\&=\frac{1}{n(n-1)}\sum_{k,\ell}
\left(\theta_i^k-\theta_i^\ell\right)\left(\theta_j^k-\theta_j^\ell\right)
\left(X_k-X_\ell\right)^2.\end{split}\end{equation*}
Expanding this expression and making use of the facts that 
$\inprod{\theta_i}{\theta_j}=0$ and $\sum_r\theta_i^r=\sum_r\theta_j^r=0$
gives that the right-hand side is equal to
$$\frac{2}{n(n-1)}\left[n\sum_k\theta_i^k
\theta_j^kX_k^2-2\left(\sum_k\theta_i^k\theta_j^kX_k\right)\left(\sum_\ell
X_\ell\right)+2S_n^iS_n^j\right]=:E_{ij}.$$
As in the case of $E_{ii}$, to estimate $\E|E_{ij}|,$ apply the triangle 
inequality to the expression above and estimate each term separately.
First,
\begin{equation*}\begin{split}
\E\left|\sum_r \theta_i^r\theta_j^rX_r^2\right|&\le
\sqrt{\E\sum_{r,s}\theta_i^r\theta_j^r\theta_i^s\theta_j^sX_r^2X_s^2}\\
&=\sqrt{\left(\E X_1^4-\E X_1^2X_2^2\right)\sum_{r}(\theta_i^r)^2(
\theta_j^r)^2}\\&\le\|\theta_i\|_4\|\theta_j\|_4\sqrt{\E X_1^4},
\end{split}\end{equation*}
where the second line follows from exchangeability and 
the fact that $\inprod{\theta_i}{\theta_j}=0$.
By the normalization, $\big|\E S_n^iS_n^j\big|\le1,$ and it remains to 
estimate the middle term.  As before, this is done by applying the 
Cauchy-Schwarz inequality, breaking up the sum by equality structure of
the indices, and using exchangeability and the orthonormality conditions
on the $\theta_i$ to simplify the result.  This
process yields
\begin{equation*}\begin{split}
\E&\left|\sum_{r,s}\theta_i^r\theta_j^rX_rX_s\right|\\&\qquad\le
\sqrt{\sum_r(\theta_i^r)^2(\theta_j^r)^2\Big[n^2\big|E X_1X_2X_3X_4\big|
+n^2\big|\E X_1^2X_2X_3\big|
+2n\big|\E X_1^3X_2\big|+n\E X_1^2X_2^2+\E X_1^4\Big]}\\&\qquad\le
\|\theta_i\|_4\|\theta_j\|_4cn\sqrt{\E X_1^4},
\end{split}\end{equation*}
for some constant $c$.  It follows that there is another constant $a$
such that for all $i\neq j$,
$$\E|E_{ij}|\le\frac{a}{n} \|\theta_i\|_4\|\theta_j\|_4\sqrt{\E X_1^4}.$$
It now follows that 
$$\frac{1}{\lambda}\sum_{i,j=1}^k\E|E_{ij}|\le ck\sqrt{\big|\E X_1X_2X_3X_4
\big|}+c'k\sqrt{\big|\E(X_1^2-1)(X_2^2-1)\big|}+c''\sqrt{\E X_1^4}
\left(\sum_i\|\theta_i\|_4\right)^2.$$
To complete the application of Theorem \ref{discrete}, it remains to estimate
$\E\big|(S_n')^i-S_n^i\big|^3.$  By \eqref{diff},
\begin{equation*}\begin{split}
\E\big|(S_n')^i-S_n^i\big|^3&=\E\left|(\theta_i^I-\theta_i^J)(X_J-X_I)\right|^3
\\&=\frac{1}{n(n-1)}\sum_{r,s}\E\left|(\theta_i^r-\theta_i^s)(X_s-X_r)\right|^3
\\&\le\frac{8\E|X_1|^3}{n(n-1)}\sum_{r,s}\left|\theta_i^r-\theta_i^s\right|^3\\
&\le\frac{8\E|X_1|^3}{n(n-1)}\sum_{r,s}\left(|\theta_i^r|+|\theta_i^s|\right)^3
\\&=\frac{8\E|X_1|^3}{n(n-1)}\Big[2n\|\theta_i\|_3^3+6\|\theta_i\|_1\Big].
\end{split}\end{equation*}
Here the third line follows from the $L_3$ triangle inequality and 
exchangeability, and the last line by expanding the cube and using the
normalization condition on $\theta_i$.
Note that, by H\"older's inequality and the fact that $\|\theta_i\|=1$ for
each $i$, 
$$\|\theta_i\|_1\le\sqrt{n}$$
and $$\|\theta_i\|_3^3\ge\frac{1}{\sqrt{n}},$$
thus the second term above can be absorbed into the first with a change in
constant.  It follows that
$$\frac{1}{\lambda}\sum_i\E\big|(S_n')^i-S_n^i\big|^3\le c\E|X_1|^3\sum_{i=1}^k
\|\theta_i\|_3^3.$$
This accounts for the remaining error term
in Theorem \ref{exch-thm}.

\end{proof}

\bibliographystyle{plain}
\bibliography{altclt.bib}

\end{document}